\documentclass[a4paper,10pt,leqno]{article}
\usepackage{amsfonts,amssymb,amsmath,amsthm}
\usepackage{fullpage}
\usepackage{epic,eepic}
\usepackage{mathrsfs}
\newtheorem{theorem}{Theorem}[section]
\newtheorem{proposition}[theorem]{Proposition}
\newtheorem{corollary}[theorem]{Corollary}
\newtheorem{lemma}[theorem]{Lemma}

\theoremstyle{remark}
\newtheorem{remark}[theorem]{{\bf Remark}}

\numberwithin{equation}{section}

\newcommand{\A}{\mathbf A}
\newcommand{\N}{{\mathbb N}}
\newcommand{\R}{{\mathbb R}}

\newcommand{\K}{{\mathcal K}}
\newcommand{\eps}{{\varepsilon}}
\newcommand{\vark}{{\varkappa}}

\renewcommand\b{\beta}
\newcommand\g{\gamma}

\renewcommand\d{\delta}
\renewcommand\l{\lambda}
\newcommand\n{\nabla}
\renewcommand\O{\Omega}

\def\var{\varphi}

\def\vark{\varkappa}

\newcommand{\esssup}{\mathop{\mathrm{ess\,sup}}}

\newcommand{\indi}{\mbox{\rm 1\hspace{-2pt}\rule[0mm]{0.2mm}{7.9pt}\hspace{ 1pt}}}
\newcommand{\ind}{\hbox{\rm 1\hskip -4.5pt 1}}

\def\part{\partial}

\def\pr{^{\prime}}

\def\XXint#1#2#3{{\setbox0=\hbox{$#1{#2#3}{\int}$}
\vcenter{\hbox{$#2#3$}}\kern-.5\wd0}}

\def\W{\rlap{$\buildrel \circ \over W$}\phantom{W}}

\begin{document}

\title{
Pointwise estimates for solutions of singular quasi-linear parabolic
equations }

\author{
{\large Vitali Liskevich}\footnote{Corresponding author}\\
\small Department of Mathematics\\
\small University of Swansea\\
\small Swansea SA2 8PP, UK\\
{\tt v.a.liskevich@swansea.ac.uk}\\
\and
{\large Igor I.\,Skrypnik}\\
\small Institute of Applied\\
\small Mathematics and Mechanics\\
\small Donetsk 83114,  Ukraine\\
{\tt iskrypnik@iamm.donbass.com}
}

\date{}

\maketitle

\setlength{\unitlength}{0.0004in}
\begingroup\makeatletter\ifx\SetFigFont\undefined%
\gdef\SetFigFont#1#2#3#4#5{%
  \reset@font\fontsize{#1}{#2pt}%
  \fontfamily{#3}\fontseries{#4}\fontshape{#5}%
  \selectfont}%
\fi\endgroup%
\renewcommand{\dashlinestretch}{30}

\begin{abstract}
For a class of singular divergence type quasi-linear parabolic  equations
with a Radon measure on the right hand side
we derive pointwise estimates for solutions via the nonlinear Wolff potentials.
\end{abstract}

\bigskip

\section{Introduction and main results}

In this note we give a parabolic version of a by now classical
result by Kilpel\"ainen-Mal\'y~\cite{KiMa}, who proved
pointwise estimates for solutions to quasi-linear $p$-Laplace type
elliptic equations with measure in the right hand side. The
estimates are expressed in terms of the nonlinear Wolff potential of
the right hand side. These estimates were subsequently extended to
fully nonlinear equations by Labutin \cite{Labutin} and fully
nonlinear and subelliptic quasi-linear equations by Trudinger and
Wang \cite{TW}. The pointwise estimates proved to be extremely
useful in various regularity and solvability problems for
quasilinear and fully nonlinear equations
\cite{KiMa,Labutin,PV1,PV2,TW}.
An immediate consequence of these estimates is the sufficient condition of local
boundedness of weak solutions which, as examples show, is optimal.

For the heat equations the
corresponding result was recently given in \cite{DM3}.
The degenerate case $p>2$ was studied recently in \cite{LSS-1}. Here
we provide the pointwise estimates for the singular supercritical
case $\frac{2n}{n+1}<p<2$.

Let $\O$ be a domain in $\R^n$, $T>0$. Let $\mu$ be a Radon measure on $\O$.
Let $\frac{2n}{n+1}<p<2$.
We are concerned with pointwise estimates for a class of non-homogeneous divergence type quasi-linear parabolic equations
of the type
\begin{equation}
\label{e0}
u_t- {\rm div}\,\A(x,t,u,\n u)=\mu\quad \text{in}\ \O_T=\O\times (0,T), \quad \O\subset \R^n.
\end{equation}
We assume that the following structure conditions are satisfied:
\begin{eqnarray}
\label{e1.3}
\A (x,t,u, \zeta)\zeta &\ge& c_1 |\zeta|^p,\quad \zeta\in \R^n,\\
\label{e1.4} |\A (x,t,u, \zeta)|&\le& c_2 |\zeta|^{p-1},
\end{eqnarray}
with some positive constants $c_1, c_2$.
%
The model example involves the parabolic $p$-Laplace equation
\begin{equation}
\label{e1.1} u_t-\Delta_pu=\mu ,\quad (x,t)\in
\O_T.
\end{equation}

\bigskip

Before formulating the main results, let us remind the reader of the definition of a weak solution
to equation \eqref{e0}.

We say that $u$ is a weak solution to \eqref{e0} if $u\in
V(\O_T):=C([0,T]; L^2_{\rm loc}(\O))\cap L^p_{\rm loc}(0,T;
W^{1,p}_{\rm loc}(\O))$ and for any compact subset $\K$ of $\O$ and
any interval $[t_1,t_2]\subset (0,T)$ the integral identity
\begin{equation}
\label{e1.12b} \int_{\K}u\var dx\Big|_{t_1}^{t_2}+\int_{t_1}^{t_2}
\int_\K \left\{-u\var_\tau+\A(x,t,u,\n u) \n\var
\right\}dx\,d\tau=\int_{t_1}^{t_2}\int_{\K} \var \mu(dx)\,d\tau.
\end{equation}
for any $\var\in W^{1,2}_{\rm loc}(0,T;C(\K))\cap L^p_{\rm loc}(0,T; \W^{1,p}(\K))$.
Note that $\var$ is required to be continuous with respect to the spatial variable
so that the right hand side of \eqref{e1.12b} is well defined.


\medskip
The crucial role in our results is played by the truncated version
of the Wolff potential defined by
\begin{equation}
\label{wolff}
W^\mu_p (x,R)=\int_0^R \left(\frac{\mu(B_r(x))}{r^{n-p}}\right)^\frac1{p-1}\frac{dr}r.
\end{equation}

In the sequel, $\g$ stands for a constant which depends only on $n,p,c_1,c_2$ which may vary from line to line.

\begin{theorem}
\label{thm1.1b}Let $u$ be a weak solution to equation~\eqref{e0}.
Let $\b=p+n(p-2)>0$.
There exists $\g>0$ depending on
$n,c_1, c_2$ and $p$,
such that
for almost all $(x_0,t_0)\in \O_T$ there exists $R_0>0$ satisfying the condition $R_0<\min\{1,t_0^\frac1\b,(T-t_0)^\frac1\b\}$
such that for all $R\le R_0$ 
the following estimates hold
\begin{enumerate}
\item[(i)]
$\displaystyle
u(x_0,t_0)\le
\g\left\{R^2+\left(\frac1{R^{p+n}}\iint_{B_R\times(0,\,T)} u_+dxdt\right)^\frac{1}{3-p}
+ W^\mu_p(x_0,2R)
\right\};$
\item[(ii)]
$\displaystyle
u(x_0,t_0)\le \g \left\{R^2+\left(\frac1{R^{n}}\esssup_{0<t<T}\int_{B_R} u_+ dx\right)
+ W^\mu_p(x_0,2R)
\right\}.$
\end{enumerate}
%
\end{theorem}

Estimate (i) is not homogeneous in $u$ which is usual for
such type of equations \cite{DiB, DiGV}. The proof of
Theorem~\ref{thm1.1b} is based on a suitable modifications of
De Giorgi's iteration technique \cite{DG} following the
adaptation of Kilpel\"ainen-Mal\'y technique \cite{KiMa} to
parabolic equations with ideas from \cite{LS2, Skr1}.

The following test for local boundedness of solutions to~\eqref{e0}
is an immediate consequence of Theorem~\ref{thm1.1b}.
\begin{corollary}
Let $u$ be a weak solution to equation~\eqref{e0}.
Let $\frac {2n}{n+1}<p<2$. Assume that there exists $R>0$ such that
\[
\sup_{x\in \O} W^\mu_p(x,R)<\infty.
\]
Then $u\in L^\infty_{\rm loc}(\O_T)$.
\end{corollary}

\begin{remark}
1. The range of $p$ in the above results is optimal for the validity of the Harnack inequality (cf.~\cite{DiB}),
however for weak solutions of the considered class it is plausible to conjecture that the results are valid
for $\frac{2n}{n+2}<p<2$, although it may require some additional global information as in \cite[Chapter~V]{DiB}.

2. Stationary solutions of \eqref{e0} solve the corresponding elliptic equation of $p$-Laplace type with measure on the right hand side,
and the Kilpel\"ainen-Mal\'y upper bound is valid for them~\cite{KiMa,MZ}. The difference with our result in Theorem~\ref{thm1.1b}(ii)
is in the additional term $R^2$ on the right hand side of the estimate. It is not clear yet whether this is a result of the employed technique
or it lies in the essence of the problem, as even for the homogeneous structure ($\mu=0$) a similar term is present in the estimates
(cf.~\cite[Chapter~V]{DiB}).
\end{remark}


\bigskip




The rest of the paper is organized as follows. In Section~2 we give auxiliary energy type estimates.
Section~3
contains the proof of Theorem~\ref{thm1.1b}. In Section~4 we provide an application
giving the global supremum estimate for the solution to a simple initial boundary problem.

\section{Integral estimates of solutions}
\label{aux}
%

\bigskip

%
%

We start with some auxiliary integral estimates for the solutions of \eqref{e0} which are formulated in the next lemma.

Define
\begin{equation}
\label{G}
G(u)=
\left\{
\begin{array}{lll}
u&\text{for}\ &\ u>1,\\
u^{2-2\l}&\text{for}\ &\ 0<u\le 1.
\end{array}
\right.
\end{equation}

Set
\[
Q_\rho^{(\d)}(y,s)=B_\rho(y)\times
(s-\d^{2-p}\rho^p,\, s+\d^{2-p}\rho^p)\subset \O_T, \quad\rho\le R.
\]
\begin{lemma}
\label{lem2.2b}
Let the conditions of Theorem~\ref{thm1.1b} be fulfilled. Let $u$ be a solution to \eqref{e0}.
Let $\xi\in C_0^\infty(Q_\rho^{(\d)}(y,s))$ be such that $\xi(x,t)=1$ for $(x,t)\in Q_{\rho/2}^{(\d)}(y,s)$ and $|\n \xi|\le \frac{2}{\rho}$,
$|\xi_t|\le c\d^{p-2}\rho^{-p}$ for some $c>0$.
Then there exists a constant $\g>0$ depending only on $n,p,c_1,c_2$ such that for any $ l,\d>0$,
any cylinder $Q_\rho^{(\d)}(y,s)$
there holds
\begin{eqnarray}
\nonumber
&&\d^2 \int_{L(t)}G\left(\frac{u(x,t)-l}{\d}\right)\xi^k(x,t)dx
+\iint_{L}
\left(1+\frac{u-l}{\d}\right)^{-1+\l}\left(\frac{u-l}{\d}\right)^{-2\l}|\n u|^p\xi^k(x,\tau)dx\,d\tau\\
\nonumber
&\le& \g \frac{\d^p}{\rho^p} \iint_L  \left(\frac{u-l}{\d}\right)\xi^{k-1}dxd\tau
+\g \frac{\d^p}{\rho^p}\iint_L \left[\left(1+\frac{u-l}{\d}\right)^{1-\l}\left(\frac{u-l}{\d}\right)^{2\l}\right]^{p-1} \xi^{k-p}dxd\tau
\\
&&\ \ \ \ \ \ \ \ \ \ \ \ \  \  +\g\, \rho^p \d^{3-p}\mu(B_\rho(y))
, \label{e2.3b}
\end{eqnarray}
where $L=Q_\rho^{(\d)}(y,s)\cap\{u>l\}$, $L(t)=L\cap \{\tau=t\}$ and $\l\in (0,1/2)$,
$k>p$.
\end{lemma}
\proof Further on, we assume that $u_t\in L^2_{\rm loc}(\O_T)$,
since otherwise we can pass to Steklov averages. First, note that
\begin{equation}
\label{e2.4b}
\int_l^u\left(1+\frac{s-l}{\d}\right)^{-1+\l}\left(\frac{s-l}{\d}\right)^{-2\l}ds\le \g \d,
\end{equation}
and
\begin{eqnarray}
\nonumber
\int_l^u dw\int_l^w\left(1+\frac{s-l}{\d}\right)^{-1+\l}\left(\frac{s-l}{\d}\right)^{-2\l}ds
\nonumber
=\int_l^u\left(1+\frac{s-l}{\d}\right)^{-1+\l}\left(\frac{s-l}{\d}\right)^{-2\l}(u-s)ds\\
\nonumber
\ge\frac12 (u-l)\int_l^{\frac{u+l}2}\left(1+\frac{s-l}{\d}\right)^{-1+\l}\left(\frac{s-l}{\d}\right)^{-2\l}ds
=\frac{\d^2}2\left(\frac{u-l}{\d}\right)\int_0^{\frac{u-l}{2\d}}(1+z)^{-1+\l}z^{-2\l}dz\\
\label{e2.5b}
\ge \g\d^2G\left(\frac{u-l}{\d}\right).
\end{eqnarray}

Let $\eta$ be the standard mollifier in $\R^n$ and as usual $\eta_\sigma(x)=\eta(\frac{x}{\sigma})$.
Set $u_\sigma = \eta_\sigma* u$. Let $\eps>0$.

Test \eqref{e1.12b} by $\var$ defined by
\begin{equation}
\var(x,t)=\left[\int_l^{u_\sigma(x,t)}\left(1+\frac{s-l}{\d}\right)^{-1+\l}\left(\eps+\frac{s-l}{\d}\right)^{-2\l}ds\right]_+\xi(x,t)^k,
\end{equation}
and $t_1=s-\d^{2-p}\rho^p$, $t_2=t$.

Using first  \eqref{e2.4b} on the right hand side, then passing to the limit $\sigma\to 0$ on the left,
after applying the  Schwarz inequality we obtain for any $t>0$

\begin{eqnarray*}
&&\int_{L(t)}\int_l^u dw\int_l^w\left(1+\frac{s-l}{\d}\right)^{-1+\l}\left(\eps+\frac{s-l}{\d}\right)^{-2\l}ds \xi^k dx\\
&+&\iint_L
\left(1+\frac{u-l}{\d}\right)^{-1+\l}\left(\eps+\frac{u-l}{\d}\right)^{-2\l}|\n u|^p\xi^k dxdt\\
&\le& \g \iint_L \int_l^udw \int_l^w\left(1+\frac{s-l}{\d}\right)^{-1+\l}\left(\eps+\frac{s-l}{\d}\right)^{-2\l}ds|\xi_t|\xi^{k-1}dxdt\\
&+&\g \frac{\d^p}{\rho^p}\iint_L \left[\left(1+\frac{u-l}{\d}\right)^{1-\l}\left(\eps+\frac{u-l}{\d}\right)^{2\l}\right]^{p-1} \xi^{k-p}dx\,dt
+\g \d^{3-p}\rho^p \mu(B_\rho(y)).
\end{eqnarray*}
Passing to the limit $\eps\to 0$ on the right hand side of the above inequality and using the Fatou lemma on the left,
by \eqref{e2.5b} we obtain the required~\eqref{e2.3b}.~\qed

\bigskip
Now set
\begin{equation}
\label{e2.6b}
\psi(x,t)=\frac1{\d}\left[\int_l^{u(x,t)}\left(1+\frac{s-l}{\d}\right)^{-\frac{1-\l}{p}}\left(\frac{s-l}{\d}\right)^{-\frac{2\l}{p}}ds\right]_+.
\end{equation}
The next lemma is a direct consequence of Lemma~\ref{lem2.2b}.
\begin{lemma}
\label{lem2.3b}
Let the conditions of Lemma~\ref{lem2.2b} be fulfilled. Let $\l<\min\{p-1,\frac{2-p}{p-1},\frac1{2}\}$.
Then
\begin{eqnarray}
\nonumber
&&\int_{L(t)}G\left(\frac{u-l}{\d}\right)\xi^k
dx+\d^{p-2}\iint_L |\n\psi|^p \xi^k dxd\tau\\
&&\le
\g \frac{\d^{p-2}}{\rho^p}\iint_L \left(1+\frac{u-l}{\d}\right)^{1-2\l(p-1)}
\left(\frac{u-l}{\d}\right)^{2\l(p-1)}\xi^{k-p}dxd\tau\ \
\label{e2.10b}
+\g \frac{\rho^p}{\d^{p-1}}\mu({B_\rho(y)}). \ \
\end{eqnarray}
\end{lemma}
\proof
Follows directly from the previous lemma and the condition on $\l$.
\qed

\section{ Proof of Theorem~\ref{thm1.1b}}
\label{bdd}

\subsection{Construction of the iteration sequences $(l_j)_j,\
(\d_j)_j, \ (\theta_{j,N}(t))_j$. }


Let $(x_0,t_0)\in \O_T$. We further assume that $R_0<\min\{1,{\rm dist}\,(x_0,\part\O),
t_0^\frac1\b,(T-t_0)^\frac1\b\}$ and $R\le R_0$. Set
\[
Q_R(x_0,t_0)=B_R(x_0)\times (t_0-R^p,t_0+R^p).
\]
Set
$\rho_j=2^{-j}R$, $B_j=B_{\rho_j}(x_0)$.

Let $\xi_j\in C_0^\infty(B_j)$, $\indi_{B_{j+1}}\le \xi_j\le
\indi_{B_j}$ with $|\n\xi_j|\le 2\rho_j^{-1}$, where here and below  $\indi_{S}$ stands for the characteristic
function of the set $S$.

The sequences of positive numbers $(l_j)_{j\in\N}$ and $(\d_j)_{j\in\N}$ are defined inductively as follows.

Set $l_{-1}=0, \ l_0=\rho_0^2,\  \d_{-1}=l_0-l_{-1}=\rho_0^2$. Let
$B\ge 1$ be a number satisfying
\begin{equation}
\label{B1}
B^{2-p}(2R_0)^\b\le
\min\{t_0,T-t_0\}
\end{equation}
which will be fixed later.

Let $j\ge 1$. Suppose we have already chosen the values $l_i, \
i=1,2,\dots, j$ and $\d_i=l_{i+1}-l_{i}$ in such a way that $\d_i\le
B\rho_i^{-n}$, $i=0,1,2,\dots,j-1$.

Let $l\in (l_j+\frac{\d_{j-1}}2,l_j+B \rho_j^{-n})$, $\d_j(l)=l-l_j$.

Take the interval
\[
I_j=\left(t_0-B^{2-p}\rho_j^\b, t_0+B^{2-p}\rho_j^\b\right).
\]
Let us divide $I_j$ into the equal parts by the points $\tau^*_{j,m}$, $m=1,2,\dots,M^*(j)$, in such a way
that
\begin{equation}
\label{j-height}
\frac14\d_{j-1}^{2-p}\rho_j^p\le \tau^*_{j,m+1} - \tau^*_{j,m} \le \frac12\d_{j-1}^{2-p}\rho_j^p, \quad \text{and}\
M^*(j)\le \g B^{2-p}\rho_j^{n(p-2)}\d_{j-1}^{p-2}.
\end{equation}

Let $\bar\theta\in C^\infty(\R)$ be such that $0\le \bar\theta \le 1$,
$\bar \theta (s)=0$ if $|s|\ge 1$, $\bar\theta(s)=1$ if $|s|\le 2^{1-p}$ and $|\bar\theta'(s)|\le \g(p)$.
Set
\[
\theta^*_{j,m}(t,l)=\bar\theta\left(\d_j(l)^{p-2}\rho_j^{-p}(t-\tau^*_{j,m})\right)
\]

Note that
\[
\sum_{m=1}^{M^*(j)}\theta^*_{j,m}(t,l)\ge \ind_{I_j}.
\]

For a fixed $j\ge 1$ and every $m=1,2,\dots M^*(j)$ we define
\begin{equation}
\label{e3.2}
\begin{split}
A^*_{j,m}(l)&=\frac{\d_j(l)^{p-2}}{\rho_j^{n+p}}\iint_{L_j}\frac{u-l_j}{\d_j(l)}\,\xi_j(x)^{k-p}
[\theta^*_{j,m}(t,l)]^{k-p} dxdt\,\\
&+\esssup_{t}\, \frac1{\rho_j^n}\int_{L_j(t)}G\left(\frac{u-l_j}{\d_j(l)}\right)\xi_j(x)^{k}
[\theta^*_{j,m}(t,l)]^{k}  dx,\\
\hbox{and set }\quad  A_j^*(l)&=\max\limits_{1\le m \le M^*(j)} A_{j,m}^*(l),
\end{split}
\end{equation}
where
\[
L_j=\{(x,t)\in \O_T\,:\, u(x,t)>l_j\}, \quad L_j(t)=\{x\in \O\,:\, u(x,t)>l_j\},
\]
$G$ is defined in \eqref{G} and $k$ will be fixed later.

It follows from \eqref{e3.2} that
\begin{equation}
\label{e3.3}
\begin{split}
A^*_{j,m}(l_j+B\rho_j^{-n})\le 3
B^{-1}\esssup_{0<t<T}\int_{B_{R_0}}|u|dx
\le 3 B^{-1}c_{R_0},
\end{split}
\end{equation}
where
\[
c_{R_0}=
\esssup_{0<t<T}\int_{B_{R_0}}|u|dx.
\]

Fix a number $\varkappa\in (0,1)$ depending on $n,p,c_1,c_2$, which will be specified later.
Choose $B$ such that
\begin{equation}
\label{B2}
 B^{-1}c_{R_0}\le \frac{\varkappa}6.
\end{equation}
This implies that $A_j^*(l_j+B\rho_j^{-n} )\le \frac{\varkappa}2$.
Note that there exists a small enough $R_0$ such that
\eqref{B1} and ~\eqref{B2} are consistent (with some $B\ge 1$).

If
\begin{equation}
\label{e3.3b}
A_j^*(l_j+\frac{\d_{j-1}}2)\le \varkappa,
\end{equation}
we set $\d_j=\frac{\d_{j-1}}2$ and $l_{j+1}=l_j+\frac{\d_{j-1}}2$.

Note that $A_j^*(l)$ is continuous as a function of $l$.
So if
\begin{equation}
\label{e3.4b}
A_j^*(l_j+\frac{\d_{j-1}}2)> \varkappa,
\end{equation}
there exists $\bar l\in (l_j+\frac{\d_{j-1}}2,l_j+B \rho_j^{-n})$
such that $A_j^*(\bar l)=\varkappa$. In this case we set
$l_{j+1}=\bar l$ and $\d_j=l_{j+1}-l_j$ in both cases.

Note that our choices guarantee that 
\begin{equation}
\label{e3.5b}
A_j(l_{j+1})\le \varkappa.
\end{equation}

Next we construct the sequence $(\theta_{j,N}(t))_{N\in\N}$ of the cut-off functions.

Since $\frac12\d_{j-1}\le \d_j\le B \rho_j^{-n}$, we can choose a subset $\{m(1),\dots,m(M(j))\}$
of the set $\{1,\dots,M^*(j)\}$ such that
\begin{eqnarray}
\label{e3.4}
\sum_{N=1}^{M(j)}\left[\theta_{j,m(N)}^*(t,l_{j+1})\right]^k \ge 1 \quad \text{for all}\
t\in (t_0-B^{2-p}\rho_j^\b, t_0+B^{2-p}\rho_j^\b),\\
\label{e3.5}
\max_{1\le N\le M(j)} A_{j,m(N)}^*(l_{j+1})\le \varkappa\quad \text{if}\ l_{j+1}=l_j+\frac{\d_{j-1}}2,\\
\label{e3.6}
\max_{1\le N\le M(j)} A_{j,m(N)}^*(l_{j+1})= \varkappa\quad \text{if}\ l_{j+1}>l_j+\frac{\d_{j-1}}2.
\end{eqnarray}

Set $\theta_{j,N}(t)=\theta_{j,m(N)}^*(t,l_{j+1})$ and
$A_{j,N}(l_{j+1})=A_{j.m(N)}^*(l_{j+1})$ for $1\le N\le M(j)$.
Then \eqref{e3.4}--\eqref{e3.6} imply that
\begin{eqnarray}
\label{e3.7}
A_j(l_{j+1})=\max_{1\le N\le M(j)}A_{j,N}(l_{j+1})\le \varkappa\quad \text{if}\ l_{j+1}=l_j+\frac{\d_{j-1}}2,\\
\label{e3.8}
A_j(l_{j+1})=\max_{1\le N\le M(j)}A_{j,N}(l_{j+1})= \varkappa\quad \text{if}\ l_{j+1}>l_j+\frac{\d_{j-1}}2,\\
\label{e3.9}
\sum_{N=1}^{M(j)}\theta_{j,m(N)}(t,l_{j+1})^k \ge 1 \quad \text{for all}\
t\in (t_0-B^{2-p}\rho_j^\b, t_0+B^{2-p}\rho_j^\b).
\end{eqnarray}
For a fixed $N$ such that $1\le N\le M(j)$ using \eqref{e3.9} define the finite sequence
$N(i),\ i=1,2,\dots,i(N,j)$,  $1\le N(i)\le M(j-1)$ so that $k-p>1$ and there exists $\g>0$ such that
\begin{equation}
\label{e3.12}
\sum_{i=1}^{i(N,j)}\theta_{j-1,N(i)}(t)^k\ge \theta_{j,N}(t),\quad i(N,j)\le \g \left(\frac{\d_j}{\d_{j-1}}\right)^{2-p}.
\end{equation}

\subsection{Main lemma}
The following lemma is a key in the Kilpel\"ainen-Mal\'y technique \cite{KiMa}.
\begin{lemma}
\label{lem3.1b}
Let the conditions of Theorem~\ref{thm1.1b} be fulfilled.
There exists $\g>0$ depending on the data, such that for all $j\ge 1$
we have
\begin{equation}
\label{e3.6b}
\d_j\le \frac12 \d_{j-1}+\g\rho_j^2
+\g \left(\frac1{\rho_j^{n-p}}\mu(B_j)
\right)^\frac1{p-1}.
\end{equation}
\end{lemma}
\proof Fix $j\ge 1$. Without loss assume that
\begin{equation}
\label{e3.7b}
\d_j>\frac12 \d_{j-1},\quad \d_j>\rho_j^2,
\end{equation}
since otherwise \eqref{e3.6b} is evident. The second inequality in \eqref{e3.7b} guarantees that $A_j(l_{j+1})=\varkappa$.

Next we claim that under conditions \eqref{e3.7b} there is a $\g>0$ such that
\begin{equation}
\label{claim}
\d_j^{p-2} \rho_j^{-(p+n)}\iint_{L_j}\xi_j \theta_{j,N}dx\,d\tau\le \g \vark.
\end{equation}
Indeed,
for $(x,t)\in L_j$ one has
\begin{equation}
\label{Lj}
\frac{u(x,t)-l_{j-1}}{\d_{j-1}}=1+\frac{u(x,t)-l_{j}}{\d_{j-1}}\ge 1.
\end{equation}
Hence
\begin{eqnarray*}
&&\d_j^{p-2} \rho_j^{-(p+n)}\iint_{L_j}\xi_j \theta_{j,N}dx\,d\tau
\le \d_j^{p-2} \rho_j^{-(p+n)}\iint_{L_j} \left(\frac{u-l_{j-1}}{\d_{j-1}}\right)\xi_{j-1}^{k}\sum_{i=1}^{i(N,j)} \theta_{j-1,N(i)}^kdx\,d\tau\\
&\le& \d_j^{p-2}\rho_j^{-p-n}i(N,j)\max_{1\le i\le M(j-1)}
\iint_{L_j}\left(\frac{u-l_{j-1}}{\d_{j-1}}\right)\xi_{j-1}^{k-p}  \theta_{j-1,N(i)}^{k-p} dx\,d\tau\\
&\le& \g\d_{j-1}^{p-2}\rho_{j-1}^{-p-n}\max_{1\le i\le M(j-1)}
\iint_{L_j}\left(\frac{u-l_{j-1}}{\d_{j-1}}\right)\xi_{j-1}^{k-p}  \theta_{j-1,N(i)}^{k-p} dx\,d\tau
\le \g \vark,
%
\end{eqnarray*}
which proves the claim.

Recall that
\begin{equation}
\label{e3.2n}
\begin{split}
A^*_{j,N}(l_{j+1})&=\frac{\d_j^{p-2}}{\rho_j^{n+p}}\iint_{\{u> l_j\}}\left(\frac{u-l_j}{\d_j}\right)
\xi_j(x)^{k-p}[\theta_{j,N}(t)]^{k-p} dxdt\,\\
&+\esssup_{t}\, \frac1{\rho_j^n}\int_{\{u> l_j\}}G\left(\frac{u-l_j}{\d_j}\right)\xi_j(x)^{k}
[\theta_{j,N}(t)]^{k}  dx.
\end{split}
\end{equation}

Let us estimate the terms in the right hand side of \eqref{e3.2n}. 
For this we decompose $L_j$ as
$L_j=L_j\pr\cup L_j^{\prime\prime}$,
\begin{equation}
\label{decomp}
L_j\pr=\left\{(x,t)\in L_j\,:\,\frac{u(x,t)-l_j}{\d_j}<\eps \right\}, \quad L_j^{\prime\prime}=L_j\setminus L\pr_j,
\end{equation}
where $\eps\in(0,1)$ depending on $n,p,c_1,c_2$ is small enough to be determined later.
By \eqref{claim} we have
\begin{equation}
\label{e3.9b}
\begin{split}
&\frac{\d_j^{p-2}}{\rho_j^{n+p}}\iint_{L\pr_j}\left(\frac{u-l_j}{\d_j}\right)\xi_j^{k-p}\theta_{j,N}^{k-p}dxd\tau\\
&\le \frac{\d_j^{p-2}}{\rho_j^{n+p}}\eps\iint_{L\pr_j}\xi_{j-1}^{k-p}\theta_{j,N} dx\,d\tau\\
&\le \eps \frac{\d_j^{p-2}}{\rho_j^{n+p}}\sum_{i=1}^{i(N,j)}\iint_{L_{j-1}}\left(\frac{u-l_{j-1}}{\d_{j-1}}\right) \xi_{j-1}^{k-p}\theta_{j-1,N(i)}dxdt\\
&\le \eps \frac{\d_j^{p-2}}{\rho_j^{n+p}}{i(N,j)}\max_{1\le i\le M(j-1)} \iint_{L_{j-1}}\left(\frac{u-l_{j-1}}{\d_{j-1}}\right) \xi_{j-1}^{k-p}\theta_{j-1,N(i)}dxdt\\
&\le \eps \left(\frac{\d_{j-1}}{\d_j}\right)^{2-p} {i(N,j)} A_{j-1}(l_j)
\le 2^{n} \eps \varkappa.
\end{split}
\end{equation}

Set
\begin{equation}
\label{psi-j}
\psi_j(x,t)=\frac1{\d_j}\left(\int_{l_j}^{u(x,t)}\left(1+\frac{s-l_j}{\d_j}\right)^{-\frac{1-\l}{p}}
\left(\frac{s-l_j}{\d_j}\right)^{-\frac{2\l}{p}}ds\right)_+,
\end{equation}
and
\[
\rho(\l)=\frac{p}{p-1-\l}.
\]

%


The following inequalities  are easy to verify
\begin{eqnarray}
\label{psi}
c\psi_j(x,t)^{\rho(\l)}\le \frac{u(x,t)-l_j}{\d_j} \ \text{for}\ (x,t)\in L_j, \ \text{and} \\
\frac{u(x,t)-l_j}{\d_j}
\le c(\eps)\psi_j(x,t)^{\rho(\l)},
\quad (x,t)\in L_j^{\prime\prime}.
\end{eqnarray}
%

Hence
\begin{eqnarray}
\frac{\d_j^{p-2}}{\rho_j^{n+p}}\iint_{L_j^{\prime\prime}}\left(\frac{u-l_j}{\d_j}\right)\xi_j^{k-p}\theta_{j,N}^{k-p}dxd\tau
\label{e3.10b}
\le
\g(\eps)\frac{\d_j^{p-2}}{\rho_j^{n+p}}\iint_{L_j^{\prime\prime}}  \psi_j^{p\frac{n+\rho(\l)}{n}}\xi_j^{k-p}\theta_{j,N}^{k-p}dxd\tau.
\end{eqnarray}
The integral in the second terms of the right hand side of \eqref{e3.10b} is estimated by using
the Gagliardo--Nirenberg inequality in the form \cite[Chapter~II,Theorem~2.2]{LaSU}
as follows
\begin{eqnarray}
\nonumber
&&\g(\eps)\frac{\d_j^{p-2}}{\rho_j^{n+p}}\iint_{L_j^{\prime\prime}}  \psi_j^{p\frac{n+\rho(\l)}{n}}\xi_j^{k-p}\theta_{j,N}^{k-p}dxd\tau\\
\label{e3.11a}
&&\le \g\,\left( \sup_t \frac1{\rho_j^n}\int_{L_j(t)}\psi_j^{\rho(\l)}\xi_j\theta_{j,N}dx\right)^\frac{p}n
\left(\frac1{\rho_j^n}\iint_{L_j}\left|\n \left(\psi_j \xi_j^\frac{(k-p)n}{p(n+\rho(\l))}\theta_{j,N}^\frac{(k-p)n}{p(n+\rho(\l))}\right)\right|^p dx\,d\tau\right).
\end{eqnarray}

Let us estimate separately the first factor in the right hand side of \eqref{e3.11a}.
\begin{eqnarray}
\nonumber
&&\sup_t \frac1{\rho_j^n}\int_{L_j(t)}\psi_j^{\rho(\l)}\xi_j\theta_{j,N}dx\ \stackrel{\text{by \eqref{psi}}}
{\le} c^{-1}\sup_t\frac1{\rho_j^n} \int_{L_j(t)}\frac{u-l_j}{\d_j}\xi_j\theta_{j,N}\,dx\\
\nonumber
&&\stackrel{\text{by \eqref{e3.7b}}}
{\le} c^{-1} \sup_t\frac1{\rho_j^n}\frac{\d_{j-1}}{\d_j}\int_{L_j(t)}\frac{u-l_{j-1}}{\d_{j-1}}\xi_j\theta_{j,N}\,dx\\
\nonumber
&&{\le} c^{-1} \sup_t\frac1{\rho_j^n}\frac{\d_{j-1}}{\d_j}\int_{L_j(t)}G\left(\frac{u-l_{j-1}}{\d_{j-1}}\right)\xi_{j-1}^k\sum_{i=1}^{i(N,j)}\theta_{j-1,N(i)}^k\,dx\\
\nonumber
&&\stackrel{\text{by \eqref{Lj}}}
\le
c^{-1} \sup_t\frac1{\rho_j^n}\frac{\d_{j-1}}{\d_j}i(N,j)\,\max_{1\le i\le M(j-1)}\int_{L_{j-1}(t)}G\left(\frac{u-l_{j-1}}{\d_{j-1}}\right)\xi_{j-1}^k \theta_{j-1,N(i)}^k\,dx\\
&&\stackrel{\text{by \eqref{e3.5b}}}
{\le} 2^n c^{-1}\left(\frac{\d_{j-1}}{\d_j}\right)^{p-1}A_{j-1}(l_j)\le \g \vark.
\label{psi1}
\end{eqnarray}

Combining \eqref{e3.10b}, \eqref{e3.11a} and \eqref{psi1} we obtain
\begin{eqnarray}
\nonumber
\frac{\d_j^{p-2}}{\rho_j^{n+p}}\iint_{L_j^{\prime\prime}}\left(\frac{u-l_j}{\d_j}\right)\xi_j^{k-p}\theta_{j,N}^{k-p}dxd\tau\\
\le
\g(\eps)\vark^\frac{p}n \, \d_j^{p-2}\rho_j^{-n}\iint_{L_j}\left|\n \left(\psi_j \xi_j^\frac{(k-p)n}{p(n+\rho(\l))}\theta_{j,N}^\frac{(k-p)n}{p(n+\rho(\l))}\right)\right|^p dx\,d\tau.
\end{eqnarray}
For the last term in the above inequality we estimate by \eqref{claim} and \eqref{psi}
\begin{eqnarray}
\nonumber
&&\d_j^{p-2}\rho_j^{-n}\iint_{L_j}\psi_j^p\left|\n  \xi_j\right|^p \theta_{j,N} dx\,d\tau
\le \g \d_j^{p-2}\rho_j^{-n-p}\iint_{L_j}\psi_j^p \xi_{j-1} \sum_{i=1}^{i(N,j)}\theta_{j-1,N(i)}^k dx\,d\tau\\
\nonumber
&& \le \g  \d_j^{p-2}\rho_j^{-n-p}i(N,j)\,\max_{1\le i\le M(j-1)}
\iint_{L_j} \left(\frac{u-l_j}{\d_j}\right)^{p-1-\l}\xi_{j-1}^{k-p}  \theta_{j-1,N(i)}^k  dx\,d\tau\\
\nonumber
&&
\stackrel{\text{by \eqref{e3.7b}}}{\le}\g  \d_{j-1}^{p-2}\rho_j^{-n-p}\max_{1\le i\le M(j-1)}\iint_{L_j} \left(\frac{u-l_{j-1}}{\d_{j-1}}\right)
\xi_{j-1}^{k-p} \theta_{j-1,N(i)}^{k-p} dx\,d\tau\\
\label{psi-p}
&&\le \g A_{j-1}(l_j) \le \g \vark.
\end{eqnarray}

By Lemma~\ref{lem2.3b}
\begin{eqnarray}
\nonumber
\frac1{\rho^n_j}\int_{L_j(t)}G\left(\frac{u-l_j}{\d_j}\right)\xi_j^k\theta_{j,N}^{k}
dx+\frac{\d_j^{p-2}}{\rho_j^n}\iint_{L_j} |\n\psi_j|^p \xi_j^k\theta_{j,N}^{k} dxd\tau\\
\nonumber
\le
\g \frac{\d_j^{p-2}}{\rho_j^{p+n}}\iint_{L_j} \left(1+\frac{u-l_j}{\d_j}\right)^{1-2\l(p-1)}
\left(\frac{u-l_j}{\d_j}\right)^{2\l(p-1)}\xi_j^{k-p}\theta_{j,N}^{k-p}dxd\tau\ \ &&\\
\label{e3.d}
+\g \frac{\rho_j^{p-n}}{\d_j^{p-1}}\mu({B_{\rho_j}(y)}). \ \ &&
\end{eqnarray}

 Using the decomposition \eqref{decomp} and the first inequality in \eqref{e3.7b} we have
\begin{eqnarray}
\nonumber
&&
\d_j^{p-2}\rho_j^{-(n+p)}\iint_{L_j}\left(1+\frac{u-l_j}{\d_j}\right)^{1-2\l(p-1)}\left(\frac{u-l_j}{\d_j}\right)^{2\l(p-1)}
\xi_j^{k-p}\theta_{j,N}^{k-p}dx\,d\tau
\\
\nonumber
&& \le \g \eps^{2\l(p-1)}\d_j^{p-2}\rho_j^{-(n+p)}\iint_{L_j}\xi_j \theta_{j,N}dx\,d\tau\\
\nonumber
&+&\g(\eps)\d_{j-1}^{p-2}\rho_{j-1}^{-(n+p)}\max_{1\le i\le M(j-1)}
\iint_{L_{j-1}}\left(\frac{u-l_{j-1}}{\d_{j-1}}\right)  \xi_{j-1}^{k-p}\theta_{j-1,N(i)}^{k-p}dx\,d\tau
 \\
 \label{A}
 && \ \ \ \ \ \ \ \ \le \g \eps^{2\l(p-1)} \vark +\g(\eps) \vark.
\end{eqnarray}
Thus we obtain the following estimate for the first term of $A_j(l_{j+1})$:
\begin{eqnarray}
\nonumber
&&\frac{\d_j^{p-2}}{\rho_j^{n+p}}\iint_{L_j}\left(\frac{u-l_j}{\d_j}\right)
\xi_j(x)^{k-p}[\theta_{j,N}(t)]^{k-p}dxd\tau\\
\label{first}
&&\le \g
\eps^{2\l(p-1)}\vark +\g(\eps)\vark^\frac pn\left(\vark +\d_j^{1-p}\rho_j^{p-n}\mu(B_j)\right).
\end{eqnarray}

Let us estimate the second term in the right hand side of \eqref{e3.2n}.
By
\eqref{e3.d} we have
\begin{eqnarray}
\nonumber
&&\sup_t \rho_j^{-n}\int_{L_j(t)} G\left(\frac{u-l_j}{\d_j}\right)\xi_j^k \theta_{j,N}^kdx\\
\nonumber
&&\le \d_j^{p-2}\rho_j^{-(n+p)}\iint_{L_j}\left(1+\frac{u-l_j}{\d_j}\right)^{1-2\l(p-1)}\left(\frac{u-l_j}{\d_j}\right)^{2\l(p-1)}\xi_j^{k-p} \theta_{j,N}^kdx\,d\tau
+ \g \d_j^{1-p}\rho_j^{p-n}\mu(B_j)\\
\nonumber
&&\text{(by using the decomposition \eqref{decomp} and \eqref{first})}\\
\label{second}
&&\le \g\eps^{2\l(p-1)}\vark +\g(\eps )\vark^\frac pn\left(\vark +\d_j^{1-p}\rho_j^{p-n}\mu(B_j)\right)
+ \g \d_j^{1-p}\rho_j^{p-n}\mu(B_j).
\end{eqnarray}

Combining \eqref{A} and \eqref{second} and choosing $\eps$ appropriately we can find $\g_1$ and $\g$ such that
\begin{equation}
\label{kappa}
\vark \le \g_1\vark^\frac pn\left(\vark +\d_j^{1-p}\rho_j^{p-n}\mu(B_j)\right) + \g \d_j^{1-p}\rho_j^{p-n}\mu(B_j).
\end{equation}
Now choosing $\vark<1$ such that
\begin{equation}
\label{kappa-fix}
\displaystyle \vark^\frac pn =\frac1{2\g_1}
\end{equation}
we have
\begin{equation}
\label{deltaj}
\d_j\le \g \left(\rho_j^{p-n}\mu(B_j)\right)^\frac1{p-1},
\end{equation}
which completes the proof of the lemma.\qed
%

In order to complete the proof of Theorem~\ref{thm1.1b} we sum up \eqref{e3.6b} with respect to $j$ from 1 to $J-1$
\begin{eqnarray}
\nonumber
&&l_J\le \g \d_0 +\g \sum_{j=1}^\infty \rho_j^2 + \g \sum_{j=1}^\infty \left(\rho_j^{p-n}\mu(B_j)\right)^\frac1{p-1}\\
\label{lJ}
&&\le \g (\d_0+R^2 +W^\mu_p(x_0,R)).
\end{eqnarray}

It remains to estimate $\d_0$. There are two cases to consider.
Either
$\d_0=\frac12\rho_0^2=\frac12 R^2$,
or
$l_1$ and $\d_0$ are defined by $A^*_0(l_1)=\varkappa$ with $\vark$ fixed in the proof of Lemma~\ref{lem3.1b} by \eqref{kappa-fix}.  It follows that there exists $m$ such that $A_{0,m}^*(l_1)=\vark$.

Using the decomposition \eqref{decomp} with $\eps$ chosen via $\vark$,
and Lemma~\ref{lem2.3b} together with \eqref{j-height} one can see that
\[
\begin{split}
 \sup_{t}R^{-n}\int_{\{u>l_0\}}G\left(\frac{u-l_0}{\d_0}\right)\xi_0(x)\theta_{0,m}(t)dx & \le \vark/2+
\g \frac{\d_0^{p-2}}{R^{n+p}}\iint_{\{u>l_0\}}\left(\frac{u-l_0}{\d_0}\right)dx\,d\tau\\
& +\g \d_0^{1-p}R^{n-p}\mu(B_R).
\end{split}
\]

Then by \eqref{e3.2n} we conclude that either
\begin{equation}
\label{kappa-1}
\frac{\d_0^{p-2}}{R^{n+p}}\iint_{\{u> l_0\}}\left(\frac{u-l_0}{\d_0}\right)\xi_0(x)^{k-p}
[\theta_{0,m}(t)]^{k-p} dxdt \ge \g,
\end{equation}
or else
\begin{equation}
\label{kappa-2}
\d_0^{p-1}\le \g R^{p-n}\mu(B_R).
\end{equation}
In case of \eqref{kappa-1} we obtain
\begin{equation}
\label{kappa-3}
\d_0^{3-p} \le \g  \iint u_+\xi_0\theta_{0,m}dxdt.
\end{equation}
Combining this with the first case
and chosing $R_0$
such that $B^{2-p}R_0^\b\le \min\{t_0,T-t_0\}$
we have
\begin{equation}
\label{e3.17b}
\d_0\le \g  \left\{\left(\frac1{R^{p+n}}\iint_{B_R\times(t_0-B^{2-p}R^\b,\, t_0+B^{2-p}R^\b)} u_+ dxdt\right)^\frac{1}{3-p}
+R^2+\left(\frac1{R^{n-p}}\mu(B_R)\right)^\frac1{p-1}\right\}.
\end{equation}

Hence the sequence $(l_j)_{j\in\N}$ is convergent, and $\d_j\to 0\,(j\to\infty)$, and we can pass to the limit $J\to\infty$ in \eqref{lJ}.
Let $l=\lim_{j\to\infty}l_j$. From \eqref{e3.5b} we conclude that
\begin{equation}
\label{e3.19b}
\frac1{\rho_j^{n+p}}\iint_{B_{\rho_j}\times(t_0-B^{2-p}\rho_j^\b,\, t_0+B^{2-p}\rho_j^\b) }(u-l)_+\le \g \vark\,  \d_j^{3-p} \to 0\quad (j\to \infty).
\end{equation}
Choosing $(x_0,t_0)$ as a Lebesgue point of the function $(u-l)_+$ we conclude that $u(x_0,t_0)\le l$ and hence $u(x_0,t_0)$
is estimated from above by
\begin{equation}
\label{estimate-1}
u(x_0,t_0)\le
\g\left\{R^2+\left(\frac1{R^{p+n}}\iint_{B_R\times(t_0-B^{2-p}R^\b,\,t_0+B^{2-p}R^\b)} u_+ dxdt\right)^\frac{1}{3-p}
+ W^\mu_p(x_0,2R)
\right\}.
\end{equation}
Applicability of the Lebesgue differentiation theorem
follows from \cite[Chap.~II,\,Sec.~3]{Guz}. This completes the proof of
the first assertion of
Theorem~\ref{thm1.1b}. Estimate (ii) is immediate consequence of \eqref{estimate-1}. \qed

%
%

\section{Example of application}
In this section we give an application of our main result, Theorem~\ref{thm1.1b}, to the
weak solution of the
following model initial
boundary value problem.
%
\begin{eqnarray}
\nonumber
u_t-\Delta_p u=\mu, \quad (x,t)\in Q=B_R\times (0,T),\\
\label{model}
u(x,t)=0, \quad (x,t)\in S=\partial B_R\times(0,T),\\
\nonumber
u(x,0)=0, \quad x\in B_R,
\end{eqnarray}
where $\mu$ is a (positive) Radon measure on $B_R$.

Before formulating the result we need to clarify what we understand
by the weak solution to the initial boundary value problem
\eqref{model}. We assume that
\[
u\in  C\big([0,T]\,;L^2(B_R)\big)\cap
L^p\big((0,T)\,;\W^{1,p}(B_R)\big)\quad \text{and}\quad
u(t,\cdot)\to 0 \ \text{in}\ L^2(B_R)\ \text{as}\ t\to 0.
\]
\begin{proposition}
Let $u$ be the weak solution to problem \eqref{model}.
If $T\ge 4^\frac1{p-1}R^\frac\b{p-1}\mu(B_R)^\frac{2-p}{p-1}$ then
\begin{equation}
\label{i-case}
 \esssup_{B_R\times(T/4,3T/4)}u(x,t)\le
\g\left\{\left(\frac{T}{R^\b}\right)^\frac1{2-p}+R^2+\sup_{B_R}W^\mu_p(x,2R)\right\}.
\end{equation}
Otherwise,
\begin{equation}
\label{ii-case} \esssup_{B_R\times(T/4,3T/4)}u(x,t)\le
\g\left\{\left(\frac{R^p}{T}\right)^\frac{n-p}\b \cdot
\left[\sup_{B_R}W^\mu_p(x,2R)\right]^\frac{p(p-1)}\b
+R^2+\sup_{B_R}W^\mu_p(x,2R) \right\}.
\end{equation}
\end{proposition}
\proof We start with a proof of the following inequality
\begin{equation}
\label{Tmu}
\sup_{0<t<T}\int_{B_R\times \{t\}} |u| dx\le T \mu(B_R).
\end{equation}
As in the proof of Lemma~\ref{lem2.2b} we need the test function for \eqref{model} to be continuous
with respect to the spatial variable to make it $\mu$-measurable.
It is clear that $u$ can be approximated by the functions
$u_n$ which are step functions with respect to $t$ on $(0,T)$ with values in $C_0^\infty(B_R)$.
Without loss $u_n\to u$ pointwise almost everywhere and in
$L^p\big((0,T)\,;\W^{1,p}(B_R)\big)$. Now, by $v_h$ we will denote the symmetric Steklov average
$v_h(x,t)=\frac1{2h}\int_{t-h}^{t+h}v(x,\tau)d\tau$.

Taking $t_1>0$ and $t_2=t\le T$ in the integral identity~\eqref{e1.12b}, testing it with
\[
\var =\left(\frac{(u_n)_h}{|(u_n)_h|+\eps}\right)_h, \quad \eps>0,
\]
and noting that $|\var|\le 1$ we obtain
\[
\int_{t_1}^t\int_{B_R}(u_h)_\tau\frac{(u_n)_h}{|(u_n)_h|+\eps}dxd\tau
+\int_{t_1}^t\int_{B_R}\big(\n u|\n u|^{p-2}\big)_h\frac{\eps\n(u_n)_h}{\big(|(u_n)_h|+\eps\big)^2} dxd\tau
\le (t-t_1)\mu(B_R).
\]
In the above inequality we first pass to the limit $n\to\infty$.
Next passing to the limit $h\to 0$ we obtain
\[
\int_{B_R}\big(|u(x,t)|-|u(x,t_1)|\big)dx-\eps\int_{B_R}\ln\frac{|u(x,t)|+\eps}{|u(x,t_1)|+\eps}dx
+\eps\int_{t_1}^t\int_{B_R}\frac{|\n u|^p}{(|u|+\eps)^2}dxd\tau\le (t-t_1)\mu(B_R).
\]
Passing to the limit $\eps\to 0$ in the above inequality and then $t_1\to 0$ we prove
\eqref{Tmu}.


Next, in the proof of Theorem~\ref{thm1.1b} we choose $B$  as $B=T\mu(B_R)$.

Let $\sigma\in (0,1]$. We divide the cylinder $Q$ into the union of
the cylinders
\[
Q_{\sigma R}(\bar x,\bar t)= B_{\sigma R}(\bar x)\times (\bar t-B^{2-p}(\sigma R)^\b, \bar t+B^{2-p}(\sigma R)^\b).
\]
We are going to apply estimate (ii) of Theorem~\ref{thm1.1b} for the points $(\bar x, \bar t)$ in the form
\begin{equation}
\label{point-via-supt}
u(\bar x,\bar t)\le \g \left\{R^2+\left(\frac1{R^{n}}\esssup_{0<t<T}\int_{B_R} u_+ dx\right)
+ W^\mu_p(\bar x,2R)
\right\}.
\end{equation}

For this $\sigma$ has to satisfy the condition
\[
B^{2-p}(\sigma R)^\b
\le 
\frac{T}4.
\]
So we choose $\sigma = \min\left\{1, \left(\frac14\right)^\frac1 \b
R^{-1} T^\frac{p-1}\b  \mu(B_R)^\frac{p-2}\b  \right\}$.

In case $T\ge 4^\frac1{p-1}R^\frac\b{p-1}\mu(B_R)^\frac{2-p}{p-1}$
we obtain \eqref{i-case} from \eqref{point-via-supt} and \eqref{Tmu}.
In the opposite case we have
\[
\begin{split}
u(\bar x,\bar t)\le \g (\sigma R)^{-n}\int_{B_{\sigma R(x_0)}}|u|dx+\g W^\mu_p(x_0,2R)+\g R^2\\
\le \g (\sigma R)^{-n} T\mu(B_R) +\g W^\mu_p(x_0,2R)+\g R^2\\
\le (\frac14)^{-\frac n{\b}} T^{1-\frac{n(p-1)}{\b}}\mu(B_R)^{1-\frac{(p-2)n}{\b}} +\g W^\mu_p(x_0,2R)+\g R^2.
\end{split}
\]
Taking supremum  over the cylinder $B_R\times (1/4T,3/4T)$ we obtain \eqref{ii-case}. \qed

\noindent {\bf Acknowledgement}.  The authors would like to
thank the anonymous referee for valuable comments and for pointing out some
discrepancies in the first version of the manuscript.

%
%
%
%
%
%
%
%
%

\begin{small}

\end{small}


\end{document}